\documentclass[12pt, a4paper, reqno]{amsart}
\usepackage{hyperref}
\hypersetup{
    colorlinks = true,
    linkcolor = [rgb]{0.1,0.2,0.5},
    citecolor = [rgb]{0.35,0.05,0.5}
}
\usepackage{%
    amsmath,
    amsthm,
    amssymb,
    tikz,
    xcolor,
    mathtools
}

\newtheorem{thm}{Theorem}[section]
\newtheorem{lem}[thm]{Lemma}
\newtheorem{prp}[thm]{Proposition}
\newtheorem{cor}[thm]{Corollary}
\theoremstyle{definition}

\theoremstyle{remark}
\newtheorem{rmk}[thm]{Remark}
\newtheorem{exm}[thm]{Example}

\numberwithin{equation}{section}

\newcommand{\CC}{\mathbb{C}}
\newcommand{\II}{\mathbb{I}}
\newcommand{\NN}{\mathbb{N}}

\newcommand{\TT}{\mathbb{T}}
\newcommand{\ZZ}{\mathbb{Z}}
\newcommand{\Aa}{\mathcal{A}}
\newcommand{\Bb}{\mathcal{B}}
\newcommand{\Cc}{\mathcal{C}}
\newcommand{\Dd}{\mathcal{D}}
\newcommand{\Ff}{\mathcal{F}}
\newcommand{\Kk}{\mathcal{K}}
\newcommand{\Ll}{\mathcal{L}}
\newcommand{\Mm}{\mathcal{M}}
\newcommand{\Oo}{\mathcal{O}}
\newcommand{\Tt}{\mathcal{T}}
\newcommand{\Uu}{\mathcal{U}}
\newcommand{\Zz}{\mathcal{Z}}

\newcommand{\clsp}{\operatorname{\overline{\text{span}}}}
\newcommand{\Aut}{\operatorname{Aut}}
\newcommand{\id}{\operatorname{id}}

\title{Lifting isomorphisms in $K$-theory through gradings of $C^*$-algebras}

\author[E.~Ruiz]{Efren Ruiz}
\address[Ruiz]{Department of Mathematics\\University of Hawaii,
Hilo\\200 W. Kawili St.\\
Hilo, Hawaii\\
96720-4091 USA}
\email[Ruiz]{ruize@hawaii.edu}

\author[A.~Sims]{Aidan Sims}
\address[Sims]{School of Mathematics and Statistics\\
University of New South Wales Sydney\\
NSW  2502\\
Australia}
\email[Sims]{aidan.sims@unsw.edu.au}

\date{\today}
\subjclass[2020]{46L05 (primary); 46L80, 46L08 (secondary)}
\keywords{Graded C*-algebras, K-theory, Cuntz--Pimsner algebra}

\thanks{This work was initiated while the authors were attending the BIRS workshop \emph{Cartan Subalgebras in Operator Algebras, and Topological Full Groups} (\#24w5175), progressed during Aidan's visit to UH Hilo; Aidan thanks Efren for his hospitality. Aidan thanks the Isaac Newton Institute for Mathematical Sciences, Cambridge, for support and hospitality during the programme \emph{Topological groupoids and their C*-algebras}, where work on this paper was undertaken. This work was supported by EPSRC grant EP/V521929/1, by a Heilbronn Distinguished Visitor Fellowship, and by ARC Discovery Project DP220101631.}

\begin{document}

\begin{abstract}
We show that every strongly $\ZZ$-graded $C^*$-algebra (equivalently, every
$C^*$-algebra carrying a strongly continuous $\TT$-action with full
spectral subspaces) is a Cuntz--Pimsner algebra, and describe subalgebras
and subspaces that can be used as the coefficient algebra and module in the
construction. We deduce that for surjective graded homomorphisms $\phi$ of
$C^*$-algebras $A$ graded by torsion-free abelian groups $H$, if the
restriction $\phi_0$ of $\phi$ to the zero-graded component $A_0$ of $A$
induces isomorphisms in $K$-theory, so does $\phi$ itself. When $H$ is free
abelian, we show how to pick out smaller subalgebras of $A_0$ on which it
suffices to check that $\phi$ induces isomorphisms in $K$-theory.
\end{abstract}

\maketitle

The purpose of this note is to record two useful general results of which
special cases have appeared in various places in the literature. The first
(Lemma~\ref{lem:grading<->CP}) is that topologically strongly $\ZZ$-graded
$C^*$-algebras are naturally realisable as Cuntz--Pimsner algebras. The
second, and the main objective (Theorem~\ref{thm:general lift K-isos}), is
that isomorphisms in $K$-theory for $C^*$-algebras lift through gradings by
discrete torsion-free abelian groups.

The first of these results has been used repeatedly: for example, in Deaconu
and Fletcher's work on iterating the Pimsner construction \cite{Dea, FlePhD,
Fle18}; in realisations of groupoid $C^*$-algebras as Cuntz--Pimsner algebras
\cite{RRS}; and in computations of $K$-theory for self-similar actions
\cite{MS2}. The second goes back at least to Elliott's computation of the
$K$-theory of noncommutative tori \cite{Ell} and has been re-proved in
various contexts for similar purposes ever since \cite{FGS, KPS5, MS2}.

The statement of our main theorem, Theorem~\ref{thm:general lift K-isos},
says that if $A$ and $B$ are strongly graded by a torsion-free abelian group
$H$ and $\phi : A \to B$ is a surjective graded homomorphism, and if the
restriction of $\phi$ to the $0$-graded subalgebra $A_0$ of $A$ induces an
isomorphism in $K$-theory, then so does $\phi$ itself. However, in many
instances it suffices to check that $\phi$ induces an isomorphism in
$K$-theory on a much smaller subalgebra of $A_0$, a phenomenon starkly
illustrated by $k$-graph $C^*$-algebras \cite[Theorem~3.7]{FGS}. So we also
formulate a more technical-looking theorem, Theorem~\ref{thm:special lift
K-isos}, which applies only to gradings by free abelian groups, but whose
hypotheses are frequently much easier to check than those of
Theorem~\ref{thm:general lift K-isos}.

Our hope is that formulating these results abstractly provides a
labour-saving device for future proofs of invariance of $K$-theory for graded
$C^*$-algebras with respect to various notions of homotopy. To illustrate how
this might go, we indicate how our results can be used in some key examples
to recover and extend known (or at least folklore) results.

\section{Background and notation}

Throughout this paper, $\II$ denotes the unit interval $[0,1]$. The
$C^*$-algebras in this paper are assumed to be $\sigma$-unital, but we do not
assume separability. Groups $H$ are, unless we say otherwise, discrete and
abelian but not necessarily countable. By a nondegenerate homomorphism of
$C^*$-algebras, we mean a homomorphism that carries approximate identities to
approximate identities (a.k.a. a $\sigma$-unital homomorphism).

\subsection{Graded \texorpdfstring{$C^*$}{C*}-algebras}

Given a discrete group $H$, an $H$-grading of a $C^*$-algebra $A$ is a
collection of closed subspaces $A_g, g \in H$ of $A$ such that $A_e$ is a
$C^*$-subalgebra of $A$, we have $A_g A_h \subseteq A_{gh}$ and $A_{g^{-1}} =
A^*_g$ for all $g,h$, we have $A = \clsp \bigcup_{g \in H} A_g$, and there is
a positive norm-decreasing linear map $E : A \to A_e$ such that $E|_{A_e}$ is
the identity and $E(A_g) = \{0\}$ for $g \not= e$. The grading is called a
\emph{strong grading} if $\clsp A_g A_h = A_{gh}$ for all $g,h$; this is
equivalent to requiring that $\clsp A_g A_{g^{-1}} = A_e$ for all $g \in H$.
The grading is called a \emph{topological grading} if the map $E$ is faithful
on positive elements.

If $H$ is amenable, then every $H$-grading is a topological grading
\cite{Exe, Rae}. If $H$ is abelian, then an $H$-grading is the same thing as
a strongly continuous action of the Pontryagin-dual group $\widehat{H}$ (the
group of characters of $H$): given an $H$-grading of $A$, the corresponding
action of $\widehat{H}$ is uniquely determined by $\gamma_\chi(a) = \chi(g)a$
for all $a \in A_g$; and given a strongly continuous action $\gamma$ of
$\widehat{H}$ on $A$, the graded subspaces of $A$ are the spectral subspaces
$A_g = \big\{\!\int_{\widehat{H}} \overline{\chi(g)}\gamma_\chi(a)\,d\chi
\colon a \in A\big\}$ and $E$ is the canonical expectation $E(a) =
\int_{\widehat{H}} \gamma_\chi(a)\,d\chi$.\footnote{More generally, gradings
by nonabelian groups correspond to coactions, and topological gradings to
normal coactions, but in this paper we only deal with gradings by abelian
groups.}

We use the following two elementary lemmas in the proof of our main theorem.

\begin{lem}\label{lem:graded subalg}
Let $G$ be a group and let $A$ be a $G$-graded $C^*$-algebra. If $H$ is a
subgroup of $G$, then $A_H \coloneqq \clsp\big(\bigcup_{h \in H} A_h\big)$ is
an $H$-graded $C^*$-subalgebra of $A$. If $A$ is strongly graded, so is $A_H$
and if $A$ is topologically graded, so is $A_H$. If $A$ is strongly graded
and $H$ is generated as a group by $F \subseteq G$, then $A_H$ is generated
as a $C^*$-algebra by $\bigcup_{h \in F} A_h$.
\end{lem}
\begin{proof}
That $A_H$ is spanned by the graded subspaces $A_h, h \in H$ is the
definition of $A_H$. That the $A_h$ are closed subspaces of $A_H$, that $A_e$
is a $C^*$-subalgebra, and that $A_h A_{h'} \subseteq A_{hh'}$ and
$A_{h^{-1}} = A_h^*$ are inherited from $A$. These properties ensure that
$A_H$ is closed under multiplication and adjoints and hence is a
$C^*$-subalgebra of $A$. The positive norm-decreasing linear map $E : A \to
A_e$ restricts to a norm-decreasing positive linear map $E_H : A_H \to A_e$.

If $A$ is strongly graded, then in particular, for $h,h' \in H$ we have $A_h
A_{h'} = A_{hh'}$ in $A$ and hence also in $A_H$. Likewise if $A$ is
topologically graded, then $E$ is faithful on positive elements, so its
restriction to $A_H$ is too, whence $A_H$ is strongly graded.

Suppose that $H$ is generated by $F$. Let $B \coloneqq C^*\big(\bigcup_{h \in
F} A_h\big)$. Then $B \subseteq A_H$ by definition. For the reverse
inclusion, fix $h \in H$. Choose $h_1 \dots h_n \in F$ and $\varepsilon_1,
\dots, \varepsilon_n \in \{1, -1\}$ such that $h = \prod^n_{i=1}
h_i^{\varepsilon_i}$. For each $i \le n$ with $\varepsilon_i = 1$ we have
$A_{h_i^{\varepsilon_i}} = A_{h_i} \subseteq B$ by definition of $B$. For $i$
such that $\varepsilon_i = -1$, we have $A_{h_i^{\varepsilon_i}} = A_{h_i}^*
\subseteq B$ because $A$ is graded and $B$ is a $C^*$-algebra. Since $A$ is
strongly graded, we have $A_h = \prod_i A_{h_i^{\varepsilon_i}}$, and it
follows that $A_h \subseteq B$. Since $h \in H$ was arbitrary, and $B$ is a
closed subspace of $A$, we deduce that $A_H \subseteq B$.
\end{proof}

\begin{lem}\label{lem:surj graded hom}
Let $G$ be a group and let $A$ and $B$ be $G$-graded $C^*$-algebras.  Suppose
$\pi \colon A \to B$ is a surjective graded homomorphism.  Then for each $g
\in G$, the restriction $\pi \colon A_g \to B_g$ is a surjection.
\end{lem}
\begin{proof}
Fix $g \in G$ and $b \in B_g$. Since $\pi$ is surjective, there exists $a \in
A$ such that $\pi(a) = b$. By \cite[Corollary~3.5]{Exe}, there is a
contractive linear map $E^A_g : A \to A_g$ that pointwise fixes $A_g$ and
annihilates $A_h$ for $h \not= g$, and there is an analogous linear
contraction $E^B_g : B \to B_g$. Since $A = \clsp\big(\bigcup_{h \in G}
A_h\big)$, continuity and linearity give $\pi \circ F^A_g = F^B_g \circ \pi$.
Hence
\[
b = F^B_g(b) = F^B_g(\pi(a)) = \pi(F^A_g(a)) \in \pi(A_g).\qedhere
\]
\end{proof}

\subsection*{Hilbert modules and \texorpdfstring{$C^*$}{C*}-correspondences}

Let $A$ be a $C^*$-algebra. By a \emph{right-Hilbert $A$-module}, we mean a
right $A$-module $X$ equipped with an $A$-valued inner product $\langle
\cdot, \cdot\rangle_A$ that is $A$-linear in the second variable, satisfies
$\langle x, y\rangle_A^* = \langle y, x\rangle_A$ and $\langle x, x\rangle_A
\ge 0$ with equality only when $x = 0$ and that is complete in the norm
$\|x\| \coloneqq \|\langle x, x \rangle_A\|^{1/2}$. We say that $X$ is
\emph{full} as a right-Hilbert $A$-module if $\clsp \langle A,A\rangle_A =
A$.

We write $\Ll(X)$ for the $C^*$-algebra
\begin{equation}\label{eq:adjointable}
\{T : X \to X \mid \text{there exists } T^* : X \to X\text{ such that } \langle T x, y\rangle_A = \langle x, T^* y\rangle_A\text{ for all }x,y\}
\end{equation}
of adjointable operators on $X$. If $T$ is adjointable, then there is a
unique operator $T^*$ as in~\eqref{eq:adjointable}, which we call the
\emph{adjoint} of $T$. For $x,y \in X$, the formula $\theta_{x,y}(z)
\coloneqq x \cdot \langle y, z\rangle_A$ defines an adjointable operator
$\theta_{x,y}$ with adjoint $\theta_{x,y}^* = \theta_{y,x}$. We write $\Kk(X)
\coloneqq \clsp\{\theta_{x,y} \colon x,y \in X\}$, which is an ideal of
$\Ll(X)$, and we call its elements \emph{compact operators} on $X$.

A \emph{$C^*$-correspondence} over $A$ is a right-Hilbert $A$-module equipped
with a left action of $A$ satisfying $\langle a \cdot x, y\rangle_A = \langle
x, a^*\cdot y\rangle_A$ for all $a$. It follows that for each $a \in A$, the
map $\phi_X(a) \colon x \mapsto a \cdot x$ is an adjointable operator, and
the resulting map $\phi_X : A \to \Ll(X)$ is a $C^*$-homomorphism. If this
homomorphism is injective we somewhat imprecisely say that $X$ has
\emph{injective left action}, and similarly if $\phi_X$ takes values in
$\Kk(X)$, we say that $X$ has \emph{compact left action}. We say that the
left action is \emph{nondegenerate} or that $X$ is \emph{essential} if $A
\cdot X$ densely spans $X$. When the module $X$ is clear from context, we
denote $\phi_X$ simply by $\phi$.

We say that a $C^*$-correspondence $X$ over $A$ is \emph{regular} if it has a
compact injective left action. Every $C^*$-algebra $A$ is a regular
$C^*$-correspondence, denoted ${_AA_A}$, over itself when endowed with
actions given by multiplication and inner product $\langle a, b\rangle_A =
a^*b$.

Given a family $X_\lambda, \lambda \in \Lambda$ of right-Hilbert $A$-modules,
we define
\[
    \bigoplus_{\lambda \in \Lambda} X_\lambda
        \coloneqq \Big\{(x_\lambda)_{\lambda \in \Lambda} \in \prod_{\lambda \in \Lambda} X_\lambda \colon \sum_{\lambda \in \Lambda} \langle x_\lambda, x_\lambda\rangle_A\text{ coverges in norm}\Big\}.
\]
This is itself a right-Hilbert module, called the \emph{direct sum} of the
$X_\lambda$, with coordinatewise right action and with inner product $\langle
(x_\lambda), (y_\lambda)\rangle_A \coloneqq \sum_\lambda \langle x_\lambda,
y_\lambda\rangle_A$. If each $X_\lambda$ is a $C^*$-correspondence, then
$\bigoplus_\lambda X_\lambda$ is also a $C^*$-correspondence under the
coordinatewise left action of $A$.

The balanced tensor product of $C^*$-correspondences $X, Y$ over $A$ is the
completion $X \otimes_A Y$ of the quotient of the algebraic tensor product $X
\odot Y$ by the subspace $\clsp\{x\cdot a \otimes y - x \otimes a\cdot y
\colon x\in X, y \in Y, a \in A\}$ in the norm coming from the inner product
such that $\langle x \otimes y, x' \otimes y'\rangle_A = \langle y, \langle
x, x'\rangle_A \cdot y'\rangle_A$. We write $X^{\otimes n}$ for the balanced
tensor product $X \otimes_A X \otimes_A \cdots \otimes_A X$ of $n$ copies
$X$, with the convention that $X^{\otimes 0} = {_AA_A}$. The \emph{Fock
space} of $X$ is the direct sum $\bigoplus_{n=0}^\infty X^{\otimes A}$.

If $A$ and $B$ are algebras, $X$ is an $A$--$A$-bimodule and $Y$ is a
$B$--$B$-bimodule, then a \emph{bimodule morphism} from $X$ to $Y$ is a pair
$(\pi_0, \pi_1)$ consisting of a homomorphism $\pi_0 : A \to B$ and a linear
map $\pi_1 : X \to Y$ such that $\pi_0(a) \cdot \pi_1(x) = \pi_1(a \cdot x)$
and $\pi_1(x)\cdot\pi_0(a) = \pi_1(x \cdot a)$ for all $a \in A$ and $x \in
X$.

If $X$ is a $C^*$-correspondence over $A$ and $Y$ is a $C^*$-correspondence
over $B$, then a \emph{morphism of correspondences} from $X$ to $Y$ is a
bimodule morphism $(\pi_0, \pi_1)$ from $X$ to $Y$ such that $\pi_0(\langle
x, x'\rangle_A) = \langle\pi_1(x), \pi_1(x')\rangle_B$ for all $x,x' \in X$.
Any such morphism induces a homomorphism $\pi^{(1)} : \Kk(X) \to \Kk(Y)$ such
that $\pi^{(1)}(\theta_{x,x'}) = \theta_{\pi_1(x), \pi_1(x')}$. If $X$ and
$Y$ are regular, then the morphism $(\pi_0, \pi_1)$ is called
\emph{covariant} if $\pi^{(1)}(\phi_X(a)) = \phi_Y(\pi_0(a))$ for all $a \in
A$.

If $X$ is a regular $C^*$-correspondence over $A$, then a
\emph{representation} of $X$ in a $C^*$-algebra $B$ is a morphism of
correspondences $(\pi, \psi) : {_AX_A} \to {_BB_B}$. The representation is
\emph{covariant} if it is covariant as a $C^*$-correspondence morphism in the
sense of the preceding paragraph. The map $\psi$ determines morphisms $\psi_n
: X^{\otimes n} \to B$ such that $\psi(x_1 \otimes_A \cdots \otimes_A x_n) =
\psi(x_1)\psi(x_2)\cdots\psi(x_n)$ for all $x_1, \dots, x_n \in X$.

\subsection*{Cuntz--Pimsner algebras}
The \emph{Toeplitz algebra} of $X$ is the universal $C^*$-algebra $\Tt_X$
generated by the image of a representation $(i_A, i_X)$ of $X$; the
\emph{Cuntz--Pimsner algebra} of $X$ is the universal $C^*$-algebra $\Oo_X$
generated by the image of a covariant representation $(j_X, j_A)$ of $X$.
There is an injective homomorphism $\iota_X : \Kk(\Ff_X) \to \Tt_X$ of the
compact operators on the Fock space of $X$ into $\Tt_X$ such that, for $x \in
X^{\otimes m}$ and $y \in X^{\otimes_n}$
\[
\iota_X(\theta_{x,y}) = \lim_i \psi_n(x)(i_A(e_i) - i_X^{(1)}(\phi_X(e_i))) \psi_m(y^*)
\]
for any approximate identity $e_i$ for $A$. The range of $\iota_X$ is an
ideal of $\Tt_X$ and is equal to the kernel of the canonical quotient map
$q_X : \Tt_X \to \Oo_X$.

The universal properties of $\Tt_X$ and $\Oo_X$ ensure that there are
strongly continuous actions $\gamma^\Tt : \TT \to \Aut(\Tt_X)$ and
$\gamma^\Oo : \TT \to \Aut(\Oo_X)$ determined by $\gamma^\Tt_z \circ i_A =
i_A$, $\gamma^\Tt_z(i_X(x)) = z i_X(x)$, $\gamma^{\Oo}_z \circ j_A = j_A$ and
$\gamma^\Oo_z(j_X(x)) = zj_X(x)$. We frequently denote both actions by
$\gamma$ and call them both the \emph{gauge action}. These determine
topological $\ZZ$-gradings on $\Tt_X$ and $\Oo_X$ as discussed above. When
the left action on $X$ is a nondegenerate compact injective left action and
$X$ is full as a right-Hilbert $A$-module, these $\ZZ$-gradings are strong
gradings.

We will need the following result of Hume, which says that the 6-term
sequence in $K$-theory induced by Pimsner's 6-term sequence in $KK$-theory
for Cuntz--Pimsner algebras is natural with respect to covariant morphisms of
correspondences (note that Hume's definition of a \emph{morphism} of
correspondences includes the covariance condition \cite[Defition~7.0.1,
Condition~(3)]{Hume}).

\begin{prp}[{\cite[Propositions 6.4.1~and~7.0.4]{Hume}}]\label{prp:Pimsner natural}
Let $A$ and $B$ be $\sigma$-unital $C^*$-algebras. Suppose that $X$ and $Y$
are full $C^*$-correspodences over $A$ and $B$ respectively with compact and
injective left actions. Suppose that $(\pi_0, \pi_1) : (A, X) \to (B, Y)$ is
a nondegenerate covariant morphism of correspondences. Let $\pi^\Oo : \Oo_X
\to \Oo_Y$ be the induced homomorphism from \cite[Theorem~3.7]{KQR}. Then
there are maps $\overline{\partial}^X_* : K_*(\Oo_X) \to K_{1-*}(A)$ and
$\overline{\partial}^Y_* : K_*(\Oo_Y) \to K_{1-*}(B)$ that make the diagram
\[
\begin{tikzpicture}[yscale=0.8]
    \node (K0Al) at (0,6) {$K_0(A)$};
    \node (K0Am) at (6,6) {$K_0(A)$};
    \node (K0OX) at (12,6) {$K_0(\Oo_X)$};
    \node (K1Ar) at (12,0) {$K_1(A)$,};
    \node (K1Am) at (6,0) {$K_1(A)$};
    \node (K1OX) at (0,0) {$K_1(\Oo_X)$};
    \draw[->] (K0Al) to node[above] {\small$\id - [X]$} (K0Am);
    \draw[->] (K0Am) to node[above] {\small$K_0(j_A)$} (K0OX);
    \draw[->] (K0OX) to node[right] {\small$\overline{\partial}^X_0$} (K1Ar);
    \draw[->] (K1Ar) to node[below] {\small$\id - [X]$} (K1Am);
    \draw[->] (K1Am) to node[below] {\small$K_1(j_A)$} (K1OX);
    \draw[->] (K1OX) to node[left] {\small$\overline{\partial}^X_1$} (K0Al);
    \node (K0Bl) at (2,4) {$K_0(B)$};
    \node (K0Bm) at (6,4) {$K_0(B)$};
    \node (K0OY) at (10,4) {$K_0(\Oo_Y)$};
    \node (K1Br) at (10,2) {$K_1(B)$};
    \node (K1Bm) at (6,2) {$K_1(B)$};
    \node (K1OY) at (2,2) {$K_1(\Oo_Y)$};
    \draw[->] (K0Bl) to node[below] {\small$\id - [Y]$} (K0Bm);
    \draw[->] (K0Bm) to node[below] {\small$K_0(j_B)$} (K0OY);
    \draw[->] (K0OY) to node[right] {\small$\overline{\partial}^Y_0$} (K1Br);
    \draw[->] (K1Br) to node[above] {\small$\id - [Y]$} (K1Bm);
    \draw[->] (K1Bm) to node[above] {\small$K_1(j_B)$} (K1OY);
    \draw[->] (K1OY) to node[left] {\small$\overline{\partial}^Y_1$} (K0Bl);
    \draw[->] (K0Al) to node[inner sep=4pt, anchor=west] {\small$K_0(\pi_0)$} (K0Bl);
    \draw[->] (K0Am) to node[left] {\small$K_0(\pi_0)$} (K0Bm);
    \draw[->] (K0OX) to node[inner sep=6pt, anchor=east] {\small$K_0(\pi^\Oo)$} (K0OY);
    \draw[->] (K1Ar) to node[inner sep=6pt, anchor=east] {\small$K_1(\pi_0)$} (K1Br);
    \draw[->] (K1Am) to node[anchor=west] {\small$K_1(\pi_0)$} (K1Bm);
    \draw[->] (K1OX) to node[inner sep=4pt, anchor=west] {\small$K_1(\pi^\Oo)$} (K1OY);
\end{tikzpicture}
\]
commute and make the six-term rectangles exact.
\end{prp}

\begin{rmk}
In Hume's result \cite[Propositions 6.4.1]{Hume}, the top-left and
bottom-right entries in the exact rectangles in the diagram of
Proposition~\ref{prp:Pimsner natural} are replaced by $K_i(J_X)$ and
$K_i(J_Y)$; our hypothesis that the left actions are compact and injective
ensures that $J_X = A$ and $J_Y = B$. Hume gives an explicit description of
the maps $\overline{\partial}^X_i$, which he denotes by
${}\mathbin{\hat{\otimes}_i} \delta_{X, PV}$ (see the paragraph immediately
preceding \cite[Propositions 6.4.1]{Hume}), but we omit the formula as we do
not need to use it. Finally, \cite[Proposition~7.0.4]{Hume} does not state
that the right-hand trapezium in the diagram of Proposition~\ref{prp:Pimsner
natural} commutes: in Hume's notation, $\delta_{PV}$ denotes only the map
${}\mathbin{\hat{\otimes}_1} \delta_{X, PV} : K_1(\Oo_X) \to K_0(J_X)$, while
the map ${}\mathbin{\hat{\otimes}_0} \delta_{X, PV} : K_0(\Oo_X) \to
K_1(J_X)$ is written $\exp_{PV}$. However, Hume's intention was that
$\delta_{PV}$ should mean either of the two maps here (we are grateful to
Hume for a helpful conversation on this point), and certainly Hume's proof
(specifically its final sentence), modified to invoke naturality of the
exponential map in the $K$-theory ideal-quotient sequence, proves that
$[\Oo_S] \mathbin{\hat{\otimes}} \exp_{PV} = \exp_{PV}
\mathbin{\hat{\otimes}} [\varphi]$ too.
\end{rmk}

\section{\texorpdfstring{$\ZZ$}{Z}-graded \texorpdfstring{$C^*$}{C*}-algebras and Cuntz--Pimsner algebras}

We show that every strongly $\ZZ$-graded $C^*$-algebra is a Cuntz--Pimsner
algebra.

\begin{lem}\label{lem:grading<->CP}
Let $A$ be a strongly $\ZZ$-graded $C^*$-algebra. For $n \in \ZZ$, let $A_n$
be the $\ZZ$-graded subspace. Let $E : A \to A_0$ be the expectation
corresponding to the grading. Then $X \coloneqq A_1$ is a full
$C^*$-correspondence over $A_0$ with compact injective left action, and there
is an isomorphism $\Psi : A \cong \Oo_{X}$ such that $\Psi|_{A_1} = j_X$ and
$\Psi|_{A_0} = j_{A_0}$.
\end{lem}
\begin{proof}
Since $A$ is graded, multiplication in $A$ defines commuting left and right
actions of $A_0$ on $X$. By definition of strong gradings, we have $A_1^* =
A_{-1}$ and $A_{-1} A_1 = A_0$. So $\langle a, b\rangle_{A_0} \coloneqq a^*
b$ is a sesquilinear, full right $A_0$-linear map from $X \times X$ to $A_0$.
It is positive definite by the $C^*$-identity in $A$. The left action of
$A_0$ on $X$ is by adjointable operators because, for $x,y \in X$ and $a \in
A_0$, we have $\langle a\cdot x, y\rangle_{A_0} = (ax)^*y = x^* (a^*y) =
\langle x, a^* \cdot y\rangle_{A_0}$. For $x,y, z \in X$ we have
$\theta_{x,y}(z) = x\cdot \langle y, z\rangle_{A_0} = x y^* z = (xy^*)\cdot z
\in A_1 A_{-1} X$. Since $A$ is strongly graded, $A_{1}A_{-1}$ is dense in
$A_0$, and it follows that the left action of $A_0$ on $A_1$ is by compacts.
It is injective because $A_0 = A_1 A_{-1}$, so for $a \in A_0 \setminus
\{0\}$ we have $0 \not= a a^* \in a A_1 A_{-1}$, forcing $a A_1 \not= \{0\}$.

By definition, the maps $\psi_0 : A_0 \to A$ and $\psi_1 : A_1 \to A$ define
a morphism of $C^*$-correspondences from ${_{A_0} X_{A_0}}$ to ${_A A_A}$.
For $x,y \in X$, the element $xy^* \in A_0$ satisfies $\psi_0(xy^*) = xy^*
\in A$. So $(\psi_0, \psi_1)$ is a representation of ${_{A_0} X_{A_0}}$ in
$A$. The homomorphism $A_0 \to \Ll(X)$ definiing the left action is the
left-multiplication map $l(a) = (x \mapsto ax)$, and so $l(xy^*)z = xy^*z =
\theta_{x,y}(z)$ for all $z$; that is $l(xy^*) = \theta_{x,y}$. So by
definition we have $\psi^{(1)}(l(xy^*)) = \psi^{(1)}(\theta_{x,y}) =
\psi_1(x)\psi_1(y)^* = xy^* = \psi_0(xy^*)$. Since $A_0$ is densely spanned
by elements of this form, it follows that $(\psi_0, \psi_1)$ is a
Cuntz--Pimsner covariant representation. It therefore induces a homomorphism
$\Psi : \Oo_{X} \to A$ such that $\Psi \circ j_X = \psi_1$ and $\Psi \circ
j_{A_0} = \psi_0$. We have $A_0 = \Psi(A_0)$; for $n > 0$ we have $A_n =
\overline{A_1^n}  = \overline{\Psi(j_X(X))^n}$, and for $n < 0$ we have $A_n
= A_{-n}^* = (\Psi(j_X(X))^n)^*$. Thus each $A_n \subseteq
\operatorname{image}(\Psi)$. Since $A = \clsp\big(\bigcup_n A_n\big)$ we
deduce that $\Psi$ is surjective. The grading on $A$ induces an action
$\beta$ of $\TT$ such that $\beta_z(a) = z^n a$ for all $a \in A_n$. Let
$\gamma$ be the gauge action on $\Oo_X$. Then for $x \in X = A_1$, we have
$\Psi(\gamma_z(j_X(x))) = \Psi (z j_X(x)) = z\psi_1(x) = zx = \beta_z(x)$;
similarly, $\Psi(\gamma_z(j_{A_0}(a))) = a = \beta_z(a)$ for $a \in A_0$.
Since $j_X(X) \cup j_{A_0}(A_0)$ generates $\Oo_X$ we deduce that $\beta_z
\circ \Psi = \Psi \circ \gamma_z$ and so the gauge-invariant uniqueness
theorem \cite[Theorem~6.2]{Kat} shows that $\Psi$ is injective.
\end{proof}

Frequently, as with graph $C^*$-algebras (see Example~\ref{ex:graphmod}
below), it is convenient to work with a smaller $C^*$-correspondence and
coefficient algebra than $A_1$ and $A_0$.

\begin{lem}\label{lem:improved grading<->CP}
Let $A$ be a strongly $\ZZ$-graded $C^*$-algebra. For $n \in \ZZ$, let $A_n$
be the $\ZZ$-graded subspace of $A$. Suppose that $\Aa \subseteq A_0$ is a
$\sigma$-unital $C^*$-subalgebra, and $X \subseteq A_1$ is a closed subspace
such that
\begin{enumerate}
    \item $X\Aa \subseteq X$ and $\clsp \Aa X = X$,
    \item $\Aa = \clsp\{x^*y \colon x,y \in X\}$,
    \item $\clsp\{xy^* \colon x,y \in X\}$ contains an approximate identity
        for $\Aa$, and
    \item $A_1 \subseteq C^*(X)$, the $C^*$-subalgebra of $A$ generated by
        $X$.
\end{enumerate}
Then, with respect to the actions by left and right multiplication and the
inner-product $\langle x, y\rangle_\Aa = x^*y$, the space $X$ is a
nondegenerate full $C^*$-correspondence over $\Aa$ with compact injective
left action, and there is a unique isomorphism $\Psi : A \cong \Oo_{X}$ such
that $\Psi|_{X} = j_X$ and $\Psi|_{\Aa} = j_{\Aa}$.
\end{lem}
\begin{proof}
By condition~(1) $X$ is a sub-$\Aa$-bimodule of $A_1$ that is nondegenerate
as a left $\Aa$-module. Condition~(2) implies that the standard $A_0$-valued
inner-product $(x,y) \mapsto x^*y$ on $A_1$ restricts to a full $\Aa$-valued
inner-product on $X$, making it a right-Hilbert $\Aa$-module. For $b \in \Aa$
and $x,y \in X$, we have $\langle bx,y\rangle_\Aa = (bx)^*y = x^*(b^*y) =
\langle x, b^* y\rangle_\Aa$, and so the left action of $\Aa$ on $X$ by
multiplication makes $X$ into a $C^*$-correspondence over $\Aa$. To see that
the left action is compact and injective, first observe that by
condition~(3), each $b \in \Aa$ satisfies $b \in \clsp\{(bx)y^* \colon xy \in
X\} \subseteq \clsp\{xy^* \colon x,y \in X\}$ by condition~(1), so $\Aa
\subseteq \clsp XX^*$. Now suppose that $bX = \{0\}$. Since $b^* \in \clsp
XX^*$, we have $bb^* \in b\clsp XX^* = \clsp (bX) X^* = \{0\}$, forcing $b =
0$ so the left action is injective. Furthermore, for $x,y, z \in X$ we have
$\theta_{x,y}(z) = x \langle y, z\rangle_\Aa = xy^* z$, so $\theta_{x,y}$ is
left-multiplication by $xy^*$. So condition~(3) shows that the left
$\Aa$-action is compact.

Let $\iota_X : X \to A_1$ and $\iota_\Aa : \Aa \to A_0$ be the inclusion
maps. The computations of the preceding paragraph show that $(\iota_X,
\iota_\Aa)$ is a covariant morphism of $C^*$-correspondences from $(X, \Aa)$
to $(A_1, A_0)$. So \cite[Theorem~3.7]{KQR} implies that there is a
homomorphism $\Theta : \Oo_X \to \Oo_{A_1}$ such that $\Theta \circ j_X =
\iota_X$ and $\Theta \circ j_\Aa = \iota_\Aa$. Condition~(4) ensures that
$j_{A_1}(A_1) \subseteq \Theta(\Oo_X)$. Lemma~\ref{lem:grading<->CP} allows
us to identify $\Oo_{A_1}$ with $A$, which is generated as a $C^*$-algebra by
$A_1$ because it is strongly graded, so $\Theta$ is surjective. Since
$\Theta$ intertwines the gauge actions on $\Oo_X$ and $\Oo_A$ and is
injective on $j_\Aa(\Aa)$ because the identification of $\Oo_{A_1}$ with $A$
carries $\iota_\Aa$ to the inclusion map $\Aa \hookrightarrow A$, the
gauge-invariant uniqueness theorem \cite[Theorem~6.2]{Kat} implies that
$\Theta$ is injective. Now $\Psi \coloneqq \Theta^{-1}$ is the desired
isomorphism.
\end{proof}

\begin{exm}\label{ex:graphmod}
To see why Lemma~\ref{lem:improved grading<->CP} is useful, consider a finite
directed graph $E$ with no sinks or sources, and the strong grading of
$C^*(E)$ induced by the gauge action. Then the coefficient algebra and module
in Lemma~\ref{lem:grading<->CP} are $A_0 = C^*(E)^\gamma$, which is the whole
AF core $\clsp\{s_\mu s^*_\nu \colon \mu, \nu \in E^*, |\mu|= |\nu|\}$, and
$A_1 = C^*(E)_1$, which is the infinite-dimensional subspace $\clsp\{s_\mu
s^*_\nu \colon \mu, \nu \in E^*, |\nu| - |\mu| = 1\}$. However,
Lemma~\ref{lem:improved grading<->CP} allows us instead to use $\Aa =
\clsp\{p_v \colon v \in E^0\} \cong C(E^0)$ and $X = \clsp\{s_e \colon e \in
E^1\} \cong C(E^1)$, the usual graph bimodule.
\end{exm}

\section{Lifting isomorphisms in \texorpdfstring{$K$}{K}-theory}

We prove that for $C^*$-algebras that are strongly graded by discrete
torsion-free abelian groups, graded homomorphisms inducing isomorphisms of
$K$-groups of $0$-graded $C^*$-subalgebras also induce isomrphisms of
$K$-groups of the enveloping $C^*$-algebras.

We present two theorems. The first, and main, theorem is
Theorem~\ref{thm:general lift K-isos}. It applies to gradings by arbitrary
torsion-free abelian groups. It has the advantage of very broad applicability
and having an uncomplicated statement, but the drawback is that the
hypothesis that $\pi|_{C_0}$ induces isomorphisms in $K$-theory can be hard
to check in examples, as highlighed by Example~\ref{ex:graphmod}.

\begin{thm}\label{thm:general lift K-isos}
Let $H$ be a discrete torsion-free abelian group and let $C$ and $D$ be
$\sigma$-unital strongly $H$-graded $C^*$-algebras. Suppose that $\pi : C \to
D$ is a surjective graded homomorphism. Suppose that $\pi_0 \coloneqq
\pi|_{C_0} : C_0 \to D_0$ induces isomorphisms $K_*(\pi_0) : K_*(C_0) \to
K_*(D_0)$. Then $K_*(\pi) : K_*(C) \to K_*(D)$ is an isomorphism.
\end{thm}

The second theorem, Theorem~\ref{thm:special lift K-isos}, which we use to
prove Theorem~\ref{thm:general lift K-isos}, applies only to gradings by free
abelian groups, and takes advantage of Lemma~\ref{lem:improved grading<->CP}
to replace $C_0$ with a subalgebra $\Cc \subseteq C_0$ on which one needs to
check that $\pi$ induces isomorphisms in $K$-theory. It has more-complicated
hypotheses, but in practice they are often relatively easy to check in the
sense that there are natural candidates for the requisite subspaces $X_g$.

For the following theorem, recall that by Lemma~\ref{lem:graded subalg}, if
$C$ is a topologically strongly $G$-graded $C^*$-algebra, and $g \in G$, then
$C_{\ZZ g} \coloneqq \clsp\big(\bigcup_{n \in \ZZ} C_{ng}\big)$ is a
(topologically) strongly $\ZZ$-graded $C^*$-subalgebra with respect to the
grading $(C_{\ZZ g})_n = C_{ng}$, and is equal to $C^*(C_g)$.

\begin{thm}\label{thm:special lift K-isos}
Let $\Lambda$ be a set, and let $\ZZ^\Lambda$ denote the free abelian group
$\bigoplus_{\lambda \in \Lambda} \ZZ$, with free abelian generators
$\{e_\lambda \colon \lambda \in \Lambda\}$. Let $C$ and $D$ be
$\sigma$-unital strongly $\ZZ^\Lambda$-graded $C^*$-algebras. Suppose that
$\pi : C \to D$ is a surjective graded homomorphism. Suppose that $\Cc
\subseteq C_0$ is a $\sigma$-unital $C^*$-subalgebra and $X_\lambda, \lambda
\in \Lambda$ a collection of closed subspaces $X_\lambda \subseteq
C_{e_\lambda}$ such that
\begin{enumerate}
    \item for each $\lambda \in \Lambda$, the pair $\Cc, X_\lambda$
        satisfies conditions (1)--(4) of Lemma~\ref{lem:improved
        grading<->CP} with respect to the enveloping $\ZZ$-graded
        $C^*$-algebra $C^\lambda \coloneqq C_{\ZZ e_\lambda} \subseteq C$,
    \item for all $\lambda,\mu \in \Lambda$  we have $\clsp X_\lambda X_\mu
        = \clsp X_\mu X_\lambda$, and
    \item $C$ is generated as a $C^*$-algebra by $\Cc \cup \bigcup_{\lambda
        \in \Lambda} X_\lambda$.
\end{enumerate}
If $\pi|_\Cc$ induces isomorphisms $K_*(\pi|_\Cc) : K_*(\Cc) \to
K_*(\pi(\Cc))$, then $\pi$ induces isomorphisms $K_*(\pi) : K_*(C) \to
K_*(D)$.
\end{thm}

\begin{rmk}\label{rmk:canonical choice}
Note that $\Cc \coloneqq C_0$ and $X_\lambda \coloneqq C_{e_\lambda}$ for all
$\lambda \in \Lambda$ always satisfy conditions (1)--(3) of
Theorem~\ref{thm:special lift K-isos}, so the theorem can always be applied
to see that if $\pi|_{C_0} : C_0 \to D_0$ induces isomorphisms in $K$-theory,
then so does $\pi : C \to D$. The point of introducing $\Cc$, $\Lambda$ and
the $X_\lambda$ is that $\Cc$ may be significantly smaller than $C_0$, and so
it may be much easier to verify that $\pi|_\Cc$ induces isomorphisms in
$K$-theory than that $\pi|_{C_0}$ does---see Example~\ref{ex:graphmod}.
 \end{rmk}

The key to the proof of Theorem~\ref{thm:special lift K-isos} is that we can
lift an isomorphism in $K$-theory through a single $\ZZ$-grading. We set this
aside as a technical lemma.

\begin{lem}\label{lem:inductive step}
Let $A$ and $B$ be $\sigma$-unital strongly $\ZZ$-graded $C^*$-algebras.
Suppose that $\pi : A \to B$ is a surjective graded homomorphism. Suppose
that $\Aa \subseteq A_0$ is a $C^*$-subalgebra and $X \subseteq A_1$ is a
closed subspace satisfying (1)--(4) of Lemma~\ref{lem:improved grading<->CP}.
Suppose that the homomorphisms $K_*(\pi|_\Aa) : K_*(\Aa) \to K_*(\pi(\Aa))$
are isomorphisms. Then $K_*(\pi) : K_*(A) \to K_*(B)$ is an isomorphism.
\end{lem}
\begin{proof}
Let $\Bb \coloneqq \pi(\Aa) \subseteq B_0$ and let $Y \coloneqq
\overline{\pi(X)} \subseteq B_1$\label{page:fnpage}\footnote{In fact,
$\pi(X)$ is already closed for the same reason that ranges of
$C^*$-homomorphisms are closed: if $\tilde{\pi} : A/\ker(\pi) \to B$ is the
induced homomorphism, then $\tilde{\pi}$ is isometric, and $\pi(X) =
\tilde{\pi}(X + (\ker(\pi) \cap X))$ is the image of the quotient Banach
space under this isometry. We include the closure symbol rather than assume
that this fact will be obvious to our readers.}. Continuity of multiplication
ensures that $\Bb$ and $Y$ satisfy conditions (1)--(3) of
Lemma~\ref{lem:improved grading<->CP}, and the assumption that $\pi$ is
surjective guarantees~(4) (see Lemma~\ref{lem:surj graded hom}).

Let $\pi_\Aa : \Aa \to \Bb$ and $\pi_X : X \to Y$ be the restrictions of
$\pi$. Then $(\pi_\Aa, \pi_X)$ is a morphism of $C^*$-correspondences because
the operations in $X$ and $Y$ are given by multiplication in $A$ and $B$, and
$\pi$ is a homomorphism. It is nondegenerate because $A_0$ and $B_0$ contain
approximate units for $A$ and $B$. We claim that it is covariant. For this,
fix $x,y,z \in X$. Then $\theta_{x,y}(z) = x \cdot \langle y, z\rangle_{\Aa}
= xy^*z$. Further fix $a \in \Aa$. By condition~(3), there are elements
$x_{n,i}, y_{n,i} \in X$ such that $a = \lim_{n \to \infty} \sum_{i=1}^{N_n}
x_{n,i} y^*_{n,i}$. For $z \in X$ we have
\[
\phi_X(a)z
    = az
    = \lim_n \sum_i x_{n,i} y^*_{n,i} z
    = \lim_n \sum_i \theta_{x_{n,i}, y_{n,i}}(z).
\]
Hence
\[
\pi^{(1)}(\phi_X(a))\pi_X(z)
    = \lim_n \sum_i \pi(x_{n,i})\pi(y_{n,i})^* \pi_X(z)
    = \pi\Big(\lim_n \sum_i x_{n,i} y_{n,i}^*\Big)\pi_X(z)
    = \pi_\Aa(a)\pi_X(z).
\]
In particular, $\big(\pi^{(1)}(\phi_X(a)) - \pi_\Aa(a)\big) \pi(zw^*) = 0$
for all $z,w \in X$. Since $\pi$ is $\sigma$-unital, condition~(3) gives
$\pi^{(1)}(\phi_X(a)) = \pi_\Aa(a)$. So $(\pi_\Aa, \pi_X)$ is covariant.

Let $\Phi : A \to \Oo_X$ and $\Psi : B \to \Oo_Y$ be the isomorphisms of
Lemma~\ref{lem:grading<->CP}. Let $\pi^{\Oo} : \Oo_X \to \Oo_Y$ be the
homomorphism induced by the correspondence morphism $(\pi_\Aa, \pi_X)$ of
\cite[Theorem~3.7]{KQR}. For $a \in \Aa$ we have $\pi^{\Oo}(\Phi(a)) =
\pi^{\Oo}(j_{\Aa}(a)) = j_{\Bb}(\pi_\Aa(a)) = \Psi(\pi(a))$ and for $x \in X$
we have $\pi^{\Oo}(\Phi(x)) = \pi^{\Oo}(j_{X}(x)) = j_{Y}(\pi_X(x)) =
\Psi(\pi(x))$. That is, $\pi^{\Oo} \circ \Phi =  \Psi \circ \pi$. This and
the commuting diagram of Proposition~\ref{prp:Pimsner natural} imply that the
diagram
\[
\begin{tikzpicture}[yscale=0.8]
    \node (K0Al) at (0,6) {$K_0(\Aa)$};
    \node (K0Am) at (6,6) {$K_0(\Aa)$};
    \node (K0OX) at (12,6) {$K_0(A)$};
    \node (K1Ar) at (12,0) {$K_1(\Aa)$};
    \node (K1Am) at (6,0) {$K_1(\Aa)$};
    \node (K1OX) at (0,0) {$K_1(A)$};
    \draw[->] (K0Al) to node[above] {\small$\id - [X]$} (K0Am);
    \draw[->] (K0Am) to node[above] {\small$K_0(\Phi^{-1} \circ j_\Aa)$} (K0OX);
    \draw[->] (K0OX) to node[right] {\small$\overline{\partial}^X_0 \circ K_0(\Phi)$} (K1Ar);
    \draw[->] (K1Ar) to node[below] {\small$\id - [X]$} (K1Am);
    \draw[->] (K1Am) to node[below] {\small$K_1(\Phi^{-1} \circ j_\Aa)$} (K1OX);
    \draw[->] (K1OX) to node[left] {\small$\overline{\partial}^X_1 \circ K_1(\Phi)$} (K0Al);
    \node (K0Bl) at (2,4) {$K_0(\Bb)$};
    \node (K0Bm) at (6,4) {$K_0(\Bb)$};
    \node (K0OY) at (10,4) {$K_0(B)$};
    \node (K1Br) at (10,2) {$K_1(\Bb)$};
    \node (K1Bm) at (6,2) {$K_1(\Bb)$};
    \node (K1OY) at (2,2) {$K_1(B)$};
    \draw[->] (K0Bl) to node[above] {\small$\id - [Y]$} (K0Bm);
    \draw[->] (K0Bm) to node[above] {\small$K_0(\Psi^{-1} \circ j_\Bb)$} (K0OY);
    \draw[->] (K0OY) to node[left] {\small$\overline{\partial}^Y_0 \circ K_0(\Psi)$} (K1Br);
    \draw[->] (K1Br) to node[below] {\small$\id - [Y]$} (K1Bm);
    \draw[->] (K1Bm) to node[below] {\small$K_1(\Psi^{-1} \circ j_\Bb)$} (K1OY);
    \draw[->] (K1OY) to node[right] {\small$\overline{\partial}^Y_1 \circ K_1(\Psi)$} (K0Bl);
    \draw[->] (K0Al) to node[inner sep=4pt, anchor=west] {\small$K_0(\pi_\Aa)$} (K0Bl);
    \draw[->] (K0Am) to node[left] {\small$K_0(\pi_\Aa)$} (K0Bm);
    \draw[->] (K0OX) to node[inner sep=6pt, anchor=east] {\small$K_0(\pi)$} (K0OY);
    \draw[->] (K1Ar) to node[inner sep=6pt, anchor=east] {\small$K_1(\pi_\Aa)$} (K1Br);
    \draw[->] (K1Am) to node[anchor=west] {\small$K_1(\pi_\Aa)$} (K1Bm);
    \draw[->] (K1OX) to node[inner sep=4pt, anchor=west] {\small$K_1(\pi)$} (K1OY);
\end{tikzpicture}
\]
commutes, and the rectangular six-term sequences in it are exact.

By hypothesis the maps $K_*(\pi_\Aa)$ are all isomorphisms. So the Five Lemma
implies that the $K_*(\pi)$ are isomorphisms.
\end{proof}

\begin{proof}[Proof of Theorem~\ref{thm:special lift K-isos}]
For this proof, we set $\Dd \coloneqq \pi(\Cc)$ and $Y_\lambda \coloneqq
\overline{\pi(X_\lambda)}$ for $\lambda \in \Lambda$\footnote{Just as on
page~\pageref{page:fnpage}, the image $\pi(X_\lambda)$ is in fact already
closed, so the closure symbol is redundant but hopefully eliminates a
potential source of confusion.}.

For each finite $F \subseteq \Lambda$, we write $A_F$ for the
$C^*$-subalgebra of $C$ generated by $\Cc \cup \bigcup_{\lambda \in F}
X_\lambda$, and $B_F = \pi(A_F)$. We first show that the restriction of $\pi$
to each $A_F$ induces isomorphisms $K_*(\pi|_{A_F}) : K_*(A_F) \to K_*(B_F)$.
We proceed by induction on $|F|$. If $|F| = 0$ then $A_F = \Cc$ and the
statement is trivial.

Now suppose that $\pi|_{A_F}$ induces isomorphisms $K_*(\pi|_{A_F}) :
K_*(A_F) \to K_*(B_F)$ whenever $|F| < k$, and fix $F \subseteq \Gamma$ with
$|F| = k$. By enumerating $F$, we can identify it with $\{1, \dots, k\}$, and
thereby identify the subgroup of $\ZZ^\Lambda$ generated by $\{e_\lambda
\colon \lambda \in F\}$ with $\ZZ^k$. We regard $\ZZ^{k-1}$ as a subgroup of
$\ZZ^k$ in the usual way. Let $\Aa \coloneqq C^*\big(\bigcup_{i \le k-1}
X_i\big)$ and $\Bb \coloneqq C^*\big(\bigcup_{i \le k-1} Y_i\big)$. Note that
$\Cc$ is $\sigma$-unital in $C_0$, which in turn is $\sigma$-unital in $C$.
Further, $\Aa$ contains $C^*(X_1)$, and therefore contains $\Cc$ by
Condition~2 of Lemma~\ref{lem:improved grading<->CP}. Hence $\Aa$ contains an
approximate identity for $C$.

The $\ZZ^\Lambda$ gradings of $C$ and $D$ restrict to $\ZZ^{k-1}$-gradings of
$\Aa$ and $\Bb$, and $\pi : C \to D$ restricts to a graded homomorphism $\pi
: \Aa \to \Bb$, which is surjective because each $Y_i = \overline{\pi(X_i)}
\subseteq \pi(\Aa)$, and $\Bb$ is generated by the $Y_i$. We claim that $\Aa$
is strongly graded. For this, fix $g = (g_1, \dots, g_{k-1})\in \ZZ^{k-1}$.
As a notational convenience, for an integer $n > 0$ and $i \le k-1$, we write
$X_i^n$ for the set $\big\{\!\prod^n_{j=1} x_j \colon x_j \in X_i\text{ for
all $j$}\big\}$ of products (calculated in $C$) of $n$ elements of $X_i$.
Given $Z \subseteq C$ we write $Z^*$ for the set $\{z^* \colon z \in Z\}$ of
adjoints (calculated in $C$) of elements of $Z$. We have
\[
\Aa_g \supseteq
    \clsp \Aa_0 \Big(\prod_{g_i > 0} X_i^{g_i}\Big)\Big(\prod_{g_i < 0} X_i^{-g_i}\Big)^*.
\]
Hence
\[
\clsp \Aa_g \Aa_{-g}
    \supseteq \clsp \Aa_0 \Big(\prod_{g_i > 0} X_i^{g_i}\Big)
        \Big(\prod_{g_i < 0} X_i^{-g_i}\Big)^*
        \Big(\prod_{g_i < 0} X_i^{-g_i}\Big)
        \Big(\prod_{g_i > 0} X_i^{g_i}\Big)^*.
\]
If at least one $g_i$ is less than zero and $j < k-1$ is the minimal element
such that $g_j < 0$, then setting $h \coloneqq g + e_j$ so that $h_i  = g_i$
for $i \not= j$ and $h_j = g_j + 1$, we have
\[
\Big(\prod_{g_i < 0} X_i^{-g_i}\Big)^* \Big(\prod_{g_i < 0} X_i^{-g_i}\Big)
    = \Big(\prod_{h_i < 0} X_i^{-h_i}\Big)^* X_j^*X_j \Big(\prod_{h_i < 0} X_i^{-h_i}\Big).
\]
An induction on $\sum_{g_i < 0} |g_i|$ using this equality and Condition~(2)
of Lemma~\ref{lem:improved grading<->CP} shows that
\[
\clsp \Big(\prod_{g_i < 0} X_i^{-g_i}\Big)^* \Big(\prod_{g_i < 0} X_i^{-g_i}\Big)
    = \Aa.
\]
Hence Condition~(1) of Lemma~\ref{lem:improved grading<->CP} shows that
\[
\clsp \Aa_g \Aa_{-g}
    \supseteq \clsp \Aa_0 \Big(\prod_{g_i > 0} X_i^{g_i}\Big)
        \Big(\prod_{g_i > 0} X_i^{g_i}\Big)^*.
\]
Since each $X_i X_i^*$ contains an approximate identity for $\Cc$ and hence
for $A_F$, we see that $\big(\prod_{g_i > 0} X_i^{g_i}\big) \big(\prod_{g_i >
0} X_i^{g_i}\big)^*$ contains an approximate identity for $\Aa$. Hence $\clsp
\Aa_g \Aa_{-g}$ contains $\Aa_0$ as required. Similarly, $\Bb$ is strongly
graded.

By definition, the strongly $\ZZ^{k-1}$-graded $C^*$-algebra $\Aa$, with
subalgebra $\Cc \subseteq \Aa_0$ and subspaces $X_i \subseteq \Aa_{e_i}$ for
$i \le k-1$ satisfy the hypotheses of the theorem, and so the inductive
hypothesis implies that $\pi|_\Aa$ induces isomorphisms $(\pi|_\Aa)_* :
K_*(\Aa) \to K_*(\Bb)$.

Write $A \coloneqq A_F$ and $B \coloneqq B_F$. By the argument above $A$ and
$B$ are strongly $\ZZ^k$-graded $C^*$-algebras and $\pi|_A \colon A \to B$ is
a surjective graded homomorphism. In particular, for each $n \in \ZZ$, the
subspace $\clsp \bigcup_{h \in \ZZ^k \colon h_k = n} A_h$ is a closed
subspace of $A$, and $A$ is topologically strongly $\ZZ$-graded with respect
to these graded subspaces, and similarly for $B$, and $\pi|_A$ is a
surjective graded homomorphism with respect to these gradings. We have $\Aa
\subseteq A_0$ by definition of the grading of $A$, and it is a
$\sigma$-unital subalgebra because it contains $\Cc$.

Let $Z \coloneqq \clsp X_{e_k} \Aa \subseteq A_1$. We claim that $A, \Aa$ and
$Z$ satisfy the hypotheses of Lemma~\ref{lem:inductive step}. We first check
the conditions of Lemma~\ref{lem:improved grading<->CP} in turn.

\textbf{Condition~(1).} We have $Z \Aa \subseteq Z$ by definition. To see
that $\clsp \Aa Z = Z$, first note that by hypothesis, for each $i < k$ we
have $\clsp X_{e_i} X_{e_k} = \clsp X_{e_k} X_{e_i}$ and hence $X_{e_i} Z =
\clsp X_{e_i} X_{e_k} \Aa \subseteq \clsp X_{e_k} X_{e_i} \Aa \subseteq
X_{e_k} \Aa = Z$. Moreover, by Condition~(3) of Lemma~\ref{lem:improved
grading<->CP}, the space $X_{e_i} X^*_{e_i} \subseteq \Aa$ contains an
approximate identity for $\Cc$, and hence $X_{e_k} \subseteq X_{e_k} X_{e_i}
X^*_{e_i}$ by Condition~(1) of Lemma~\ref{lem:improved grading<->CP}. Thus
\[
\clsp X_{e_i}^* Z = \clsp X_{e_i}^* X_{e_k} X_{e_i} X_{e_i}^* \Aa
    = \clsp X_{e_i}^* X_{e_i} X_{e_k} X_{e_i}^* \Aa.
\]
We have $X_{e_i}^* X_{e_i} \subseteq \Cc$ by hypothesis and $X_{e_i}^* \Aa
\subseteq \Aa$ by definition of $\Aa$. Hence
\[
    \clsp X_{e_i}^* Z \subseteq \clsp \Cc X_{e_k} \Aa.
\]
Since $\clsp \Cc X_{e_k} = X_{e_k}$ by Condition~1 of Lemma~\ref{lem:improved
grading<->CP}, we obtain
\[
\clsp X_{e_i}^* Z \subseteq Z.
\]
Since $\Aa$ is spanned by products of elements of the $X_{e_i}, X_{e_i}^*$,
it follows that $\clsp \Aa Z \subseteq Z$. Since $\Aa$ contains $\Cc$ and
$\clsp \Cc X_{e_k} = X_{e_k}$ by Condition~1 of Lemma~\ref{lem:improved
grading<->CP}, the reverse containment is trivial, so $\clsp \Aa Z = Z$.

\textbf{Condition~(2).} We have $\clsp\{x^*y \colon x,y \in Z\} = \clsp \Aa
X_{e_k}^* X_{e_k} \Aa \subseteq \clsp \Aa \Cc \Aa$. This is contained in
$\Aa$ because $\Cc \subseteq \Aa$ and the reverse containment holds because
$\Cc$ contains an approximate identity for $C_0$ and hence for $C$.

\textbf{Condition~(3).} We have $\clsp\{xy^* \colon x,y \in Z\} = \clsp
X_{e_k} \Aa X_{e_k}^*$. Condition~(1) of Lemma~\ref{lem:improved
grading<->CP} implies that this contains $\clsp X_{e_k} \Cc X_{e_k}^*$, and
hence contains $X_{e_k} X^*_{e_k}$. Now Condition~(3) of
Lemma~\ref{lem:improved grading<->CP} for $X_{e_k}$ and $\Cc$ ensures that
$\clsp\{ xy^* \colon x,y \in Z\}$ contains an approximate identity for $\Cc$,
and hence for $\Aa$.

\textbf{Condition~(4).} We have already seen that $C^*(Z)$ contains $\Aa$,
which in turn contains $X_{e_i}$ for each $i < k$. It contains $X_{e_k}$
because $X_{e_k} \subseteq Z$. Since $A$ is generated by $\bigcup_{i \le k}
X_{e_i}$, we deduce that $A \subseteq C^*(Z)$, and in particular $A_1
\subseteq C^*(Z)$.

To prove the claim, it remains to check that $\pi|_{\Aa}$ induces
isomorphisms in $K$-theory, which, as we already saw, follows from the
inductive hypothesis. So $A, \Aa$ and $Z$ satisfy the hypotheses of
Lemma~\ref{lem:inductive step} as claimed.

Now Lemma~\ref{lem:inductive step} implies that $\pi|_A = \pi|_{A_F}$ induces
isomorphisms $(\pi|_{A_F})_* : K_*(A_F) \to K_*(B_F)$.

To extend this to the map $\pi : C \to D$, choose an increasing net
$(F_\alpha)_{\alpha \in I}$ of finite subsets $F_\alpha \subseteq \Lambda$
such that $\Lambda = \bigcup_{\alpha \in I} F_\alpha$. For each $\alpha$, let
$A_\alpha \coloneqq A_{F_\alpha} \subseteq C$, and $B_\alpha \coloneqq
B_{F_\alpha} \subseteq D$. Fix $\alpha \in I$. We saw above that $\pi_\alpha
\coloneqq \pi|_{A_\alpha}$ induces isomorphisms $(\pi_\alpha)_* :
K_*(A_\alpha) \to K_*(B_\alpha)$. For each $\alpha$, let $\imath_\alpha :
A_\alpha \to C$ and $\jmath_\alpha : B_\alpha \to D$ be the inclusion maps.

We claim that $K_*(\pi_\alpha)\big(\ker(K_*(\imath_\alpha))\big) =
\ker(K_*(\jmath_\alpha))$ for all $\alpha$. To see this, recall that for each
$\alpha$, we have $\ker(K_*(\imath_\alpha)) = \bigcup_{\beta \ge \alpha}
\ker(K_*(\imath_{\alpha,\beta}))$ and similarly $\ker(K_*(\imath_\alpha)) =
\bigcup_{\beta \ge \alpha} \ker(K_*(\imath_{\alpha,\beta}))$. So it suffices
to show that $K_*(\pi_\alpha)\big(\ker(K_*(\imath_{\alpha,\beta}))\big) =
\ker(K_*(\jmath_{\alpha,\beta}))$ whenever $\alpha \le \beta$ in $I$. Since
$\pi_\beta \circ \imath_{\alpha,\beta} = \jmath_{\alpha,\beta} \circ
\pi_\alpha$, we have $K_*(\pi_\beta) \circ K_*(\imath_{\alpha,\beta}) =
K_*(\jmath_{\alpha,\beta}) \circ K_*(\pi_\alpha)$. Since $K_*(\pi_\alpha)$
and $K_*(\pi_\beta)$ are isomorphisms, we obtain
$K_*(\pi_\alpha)\big(\ker(K_*(\imath_{\alpha,\beta}))\big) =
\ker(K_*(\jmath_{\alpha,\beta}))$ for all $\beta \geq \alpha$. This proves
the claim.

It follows that $K_*(\pi_\alpha)$ induces an isomorphism between the
quotients $K_*(A_\alpha)/\ker(K_*(\imath_\alpha))$ and
$K_*(B_\alpha)/\ker(K_*(\jmath_\alpha))$, which coincides with
$K_*(\pi)|_{\operatorname{Im}(K_*(\imath_\alpha))}$, from
$\operatorname{Im}(K_*(\imath_\alpha))$ to
$\operatorname{Im}(K_*(\jmath_\alpha))$. Continuity of $K$-theory shows that
$K_*(C)$ is the increasing union $\bigcup_\alpha
K_*(\imath_\alpha)(K_*(A_\alpha))$ and $K_*(D)$ is the increasing union
$\bigcup_\alpha K_*(\jmath_\alpha)(K_*(B_\alpha))$. We have just seen that
$K_*(\pi)$ restricts to isomorphisms $K_*(\imath_\alpha)(K_*(A_\alpha)) \to
K_*(\jmath_\alpha)(K_*(B_\alpha))$. Thus, it is an isomorphism $K_*(C) \to
K_*(D)$.
\end{proof}

\begin{proof}[Proof of Theorem~\ref{thm:general lift K-isos}]
Choose an increasing net $(H_{\alpha})_{\alpha \in I}$ of finitely generated
subgroups of $H$ whose union is $H$. For each $\alpha$, the fundamental
theorem of finitely generated abelian groups shows that $H_\alpha \cong
\ZZ^{k_\alpha}$ for some $k_\alpha \ge 0$; so we can choose a finite set
$\{e_1, \dots, e_{k_\alpha}\}$ of free abelian generators for $H_\alpha$. Let
$A_\alpha \coloneqq \clsp \bigcup_{g \in H_\alpha} C_g$. Then $A_\alpha$ is a
strongly $\ZZ^{k_\alpha}$-graded $C^*$-algebra, and $\pi$ restricts to a
surjective graded homomorphism $\pi_\alpha \colon A_\alpha \to B_\alpha
\coloneqq \clsp \bigcup_{g \in H_\alpha} D_g$. By Remark~\ref{rmk:canonical
choice}, the subalgebra $\Aa \coloneqq C_0$ and subspaces $X_{e_i} = C_{e_i}$
satisfy the hypotheses of Theorem~\ref{thm:special lift K-isos}, which
therefore implies that $\pi_\alpha$ induces isomorphisms $(\pi_\alpha)_* :
K_*(A_\alpha) \to K_*(B_\alpha)$. We can now proceed precisely as in the
final paragraph of the proof of Theorem~\ref{thm:special lift K-isos}.
\end{proof}

\section{Examples}

We illustrate the sorts of situations where Theorems \ref{thm:general lift
K-isos}~and~\ref{thm:special lift K-isos} could be useful through three
examples. The first two are essentially folklore results, and the third is
only a slight generalisation of \cite[Theorem~3.9]{FGS}; we do not claim
originality---we just want to indicate the utility of the theorem.

\subsection{Homotopy equivalence}\label{sec:homotopy equiv}

We continue to write $\II$ for the unit interval $\II = [0,1]$. Recall that
homomorphisms $\phi_0, \phi_1 : C \to D$ of $C^*$-algebras are said to be
\emph{homotopic} if $t \mapsto \phi_t : \{0,1\} \to \operatorname{Hom}(C, D)$
extends to a pointwise-continuous\footnote{by pointwise continuous we mean
that $t \mapsto \phi_t(a)$ is continuous for each $a \in C$} map $t \mapsto
\phi_t : \II \to \operatorname{Hom}(C, D)$. By \cite[Proposition~3.2.6]{RLL},
if $\phi_0, \phi_1 : C \to D$ are homotopic, then $K_*(\phi_0) =
K_*(\phi_1)$.

Recall also that $C^*$-algebras $C$ and $D$ are homotopic if there are
homomorphisms $\phi : C \to D$ and $\psi : D \to C$ such that $\phi \circ
\psi$ is homotopic to $\id_D$ and $\psi \circ \phi$ is homotopic to $\id_C$;
it then follows from the preceding paragraph that $K_*(\phi)$ and $K_*(\psi)$
are mutually inverse, and in particular $K_*(C) \cong K_*(D)$ via
$K_*(\phi)$. Our theorem lifts this result to a notion of homotopy
equivalence for graded $C^*$-algebras.

\begin{cor}\label{cor:homotopy}
Let $H$ be a discrete torsion-free abelian group. Suppose that $C$ and $D$
are $\sigma$-unital strongly $H$-graded $C^*$-algebras and that $\phi : C \to
D$ and $\psi : D \to C$ are graded homomorphisms. Suppose that $\phi$ is
surjective and that $\phi|_{C_0} \circ \psi_{D_0}$ is homotopic to
$\id_{D_0}$ and $\psi|_{D_0} \circ \phi|_{C_0}$ is homotopic to $\id_{C_0}$.
Then $K_*(\phi) : K_*(C) \to K_*(D)$ is an isomorphism.
\end{cor}
\begin{proof}
By \cite[Proposition~3.2.6]{RLL}, we have $K_*(\phi|_{C_0}) \circ
K_*(\psi_{D_0}) = \id_{K_*(D_0)}$ and $K_*(\psi_{D_0}) \circ K_*(\phi|_{C_0})
= \id_{K_*(C_0)}$. In particular, $\phi : C \to D$ is a surjective graded
homomorphism such that $\phi_0 \coloneqq \phi|_{C_0} : C_0 \to D_0$ induces
isomorphisms $K_*(\phi_0) : K_*(C_0) \to K_*(D_0)$. Hence
Theorem~\ref{thm:general lift K-isos} implies that the homomorphisms
$K_*(\phi) : K_*(C) \to K_*(D)$ are isomorphisms.
\end{proof}

\subsection{Twisted crossed products}\label{sec:twisted cp}

Our motivating example was Elliott's proof that all rank-$k$ noncommutative
tori have the same $K$-groups because the group of $\TT$-valued cocycles on
$\ZZ^k$ is path connected. In this subsection we show how to recover a
natural extension of this theorem from Theorem~\ref{thm:general lift K-isos}.
There are, of course, more general notions of twisted actions than the one
discussed here \cite{PR} and the argument can be extended more or less
verbatim; but we are only trying to indicate the utility of our theorem, so
have opted for a less than maximally general set-up that required relatively
little background.

A normalised \emph{2-cocycle} on a discrete group $H$ is a map $\sigma : H
\times H \to \TT$ such that
\[
    \sigma(g,h)\sigma(gh,k) = \sigma(g,hk)\sigma(h,k)\quad\text{and}
    \quad \sigma(e, g) = \sigma(g, e) = 1\quad\text{for all $g,h,k \in H$.}
\]
A \emph{path of normalised cocycles on $H$} is a map $\sigma : \II \times H
\times H \to \TT$,  $(t, g, h) \mapsto \sigma_t(g,h)$ such that $t \mapsto
\sigma_t(g,h)$ is continuous for all $g,h$ and $(g,h) \mapsto \sigma_t(g,h)$
is a normalised $2$-cocycle for each $t$.

If $H$ is a discrete abelian group, $\alpha : H \to \Aut(A)$ is an action on
a unital $C^*$-algebra, and $\sigma : H \times H \to \TT$ is a normalised
$2$-cocycle, then a covariant represenation of $(A, \alpha, \sigma)$ in a
unital $C^*$-algebra $B$ consists of a homomorphism $\pi : A \to B$ and a map
$U : H \to \Uu(B)$ such that $\pi(\alpha_g(a)) = U(g) \pi(a) U(g)^*$ and
$U(g)U(h) = \sigma(g,h) U(gh)$ for all $a \in A$ and $g,h \in G$. The
\emph{twisted crossed product} $A \rtimes_{\alpha,\sigma} H$ is the universal
$C^*$-algebra generated by a covariant representation $(i_A, i_G)$ of $(A,
\alpha, \sigma)$.

\begin{cor}\label{cor:crossedprod}
Let $H$ be a discrete countably generated torsion-free abelian group and let
$A$ be a unital $C^*$-algebra. Suppose that $\alpha : \II \times H \to
\Aut(A)$, $(t, g) \mapsto \alpha^t_g$ is a pointwise-continuous family of
maps, and $\sigma : \II \times G \times G \to \TT$, $(t, g, h) \mapsto
\sigma_t(g,h)$ is a path of normalised $2$-cocycles. Then $K_*(A
\rtimes_{\alpha^0, \sigma_0} H) \cong K_*(A \rtimes_{\alpha^1, \sigma_1} H)$.
\end{cor}
\begin{proof}
Define $\beta : H \to \Aut(C(\II, A))$ by $\beta_g(f)(t) = \alpha^t_g(f(t))$.
Let $C$ be the universal $C^*$-algebra generated by a unital homomorphism
$\pi : C(\II, A) \to C$ and a map $U : H \to \Uu(C)$ such that  $U_g \pi(f)
U^*_g = \pi(\beta_g(f))$ and $U_g U_h = \pi(t \mapsto \sigma_t(g,h)1_A)
U_{gh}$ for all $f \in C(\II, A)$ and $g,h \in H$. Identify $C(\II)$ with the
subalgebra of $C(\II, A)$ consisting of functions taking values in $\CC 1_A$;
observe that since each $\alpha^t_g$ is an automorphism, we have $\beta_g(f)
= f$ for all $f \in C(\II)$, and hence $U_g \pi(f) U_g^* = \pi(f)$ for all $f
\in C(\II)$; that is, the $U_g$ commute with elements of $C(\II)$.

For $t \in \II$, write $I_t$ for the ideal of $C$ generated by $C_0(\II
\setminus \{t\})$. Since $C_0(\II)$ is central, $I_t = \{\pi(f)c \colon c \in
C, f \in C(\II \setminus\{t\})\}$. For $a \in A$, let $a 1_\II$ denote the
constant function $t \mapsto a$ in $C(\II, A)$. Write $j^t_A : A \to C/I_t$
for the map $j_t(a) = \pi(a1_\II) + I_t$ and write $j^t_H : H \to \Uu(C/I_t)$
for the map $j^t(g) = U_h + I_t$. Then $(j^t_A, j^t_H)$ is a covariant
representation of $(A, \alpha^t, \sigma_t)$, and so induces a homomorphism
$\phi_t : A \rtimes_{\alpha^t, \sigma_t} H \to C/I_t$.

Now define $\psi : C(\II, A) \to A \rtimes_{\alpha^t, \sigma_t} H$ by
$\psi(f) = f(t)$. Then $(\psi, i_A)$ satisfies the defining relations for
$C$, and hence induces a homomorphism $(\psi \times i_A) : C \to A
\rtimes_{\alpha^t, \sigma_t} H$ such that $(\psi \times i_A) \circ \pi =
\psi$ and $(\psi \times i_A)\circ U = i_A$. In particular $C_0(\II \setminus
\{t\}) \subseteq \ker(\psi \times i_A)$, and so $(\psi \times i_A)$ descends
to a homorphism $(\psi \times i_A)_t : C/I_t \to A \rtimes_{\alpha^t,
\sigma_t} H$. Direct computation with generators shows that $\phi_t$ and
$(\psi \times i_A)_t$ are mutually inverse, and in particular $\phi_t$ is an
isomorphism. The existence of the regular representation
\cite[Definition~3.10 and Theorem~3.11]{PR} shows that $i_A : A \to A
\rtimes_{\alpha^t, \sigma_t} H$ is injective, and we deduce that $\pi :
C(\II, A) \to C$ is injective.

It therefore sufffices to show that $K_*(C_0) \cong K_*(C_1)$. For this, it
suffices to show that $K_*(C) \cong K_*(C_t)$ for each $t$. So fix $t \in
\II$. The universal property of $C$ gives rise to a strongly continuous
action $\widehat{\beta}$ of $\widehat{H}$ such that $\widehat{\beta}_\chi
(\pi(f)) = \pi(f)$ and $\widehat{\beta}_\chi(U_g) = \chi(g)U_g$ for $f \in
C(\II, A)$ and $g \in H$. Since each $\widehat{\beta}_\chi$ pointwise fixes
$C_0(\II\setminus\{t\})$, the action $\beta$ descends to an action
$\widehat{\beta}^t$ on $C/I_t$. These actions induce gradings of $C$ and
$C/I_t$ whose $g$-graded subspaces are $C_g = \pi(C(\II, A))U_g$ and
$(C/I_t)_g = \pi(C(\II, A))U_g + I_t$. In particular, the quotient map $q_t :
c \mapsto c + I_t$ from $C$ to $C/I_t$ is a surjective graded homomorphism.
Since $U_{g^{-1}} \in C_{g^{-1}}$, we have $C_0 = \pi(C(\II, A)) = \pi(C(\II,
A)) U_g U_{g^{-1}} \subseteq C_g C_{g^{-1}}$ for all $f$. Hence the grading
on $C$ is a strong grading, and taking quotients throughout shows that the
grading on $C/I_t$ is also a strong grading.

The $0$-graded component of $C$ is $C_0 = \pi(C(\II, A))$. We saw above that
$\pi$ is injective, and so $C_0 \cong C(\II, A)$, and $C_0/I_t \cong A$, and
these isomorphisms intertwine $q_t$ with the evaluation map $\varepsilon_t :
f \mapsto f(t)$ from $C(\II, A) \to A$. Since $\II$ is contractable, the maps
$K_*(\varepsilon_t) : K_*(C(\II, A)) \to K_*(A))$ are isomorphisms. So the
maps $K_*(q_t|_{C_0})$ are isomorphisms $K_*(C_0) \to K_*(C_0/I_t)$, and then
Theorem~\ref{thm:general lift K-isos} implies that the maps $K_*(q_t) :
K_*(C) \to K_*(C/I_t)$ are isomorphisms.
\end{proof}

\subsection{Homotopies of product systems}\label{sec:ps homotopy}

We show how to recover and extend \cite[Theorem~3.7]{FGS}. For detailed
background on product systems and their homotopies, we refer the reader to
\cite[Sections 2~and~3]{FGS}. Recall that $\II$ denotes the unit interval
$\II = [0,1]$.

Let $A$ be a $C^*$-algebra and fix a set $\Lambda$. We write $\NN^\Lambda$
for the canonical positive cone in $\ZZ^\Lambda = \bigoplus_{\lambda \in
\Lambda} \ZZ$, and we write $\{e_\lambda \colon \lambda \in \Lambda\}$ for
the canonical generators of $\NN^\Lambda$. By a \emph{regular product system}
over $\NN^\Lambda$ of $C^*$-correspondences over $A$, we mean a
(multiplicative) semigroup $X$ and a map $d : X \to \NN^\Lambda$ such that
$X_0 \coloneqq d^{-1}(0) = A$, for each $p \in \NN^\Lambda$ the subset $X_p =
d^{-1}(p)$ is a nondegenerate full regular $C^*$-correspondence over $A$, and
such that for all $p,q \in \NN^\Lambda$ there is an isomorphism $m_{p,q} :
X_p \otimes_{A} X_q \to X_{p+q}$ such that $m_{p,q}(x \otimes y) = xy$ for
all $x \in X_p$ and $y \in X_q$. We typically denote $X_0$ by $A$ when
thinking of it as a $C^*$-algebra, and by $X_0$ when thinking of it as the
standard correspondence ${_AA_A}$ over $A$.

A \emph{homotopy of regular product systems} over $A$ is a regular product
system $X$ over $\NN^\Lambda$ of $C^*$-correspondences over $C(\II, A)$ such
that the canonical copy of $C(\II)$ in $\Zz\Mm(C(\II, A))$ is central in $X$
in the sense that for $a, b \in C(\II, A)$ and $x \in X_p$, we have $(fa)
\cdot x \cdot b = a \cdot x \cdot (bf)$; since $X$ is regular and by
Cohen--Hewitt factorisation, elements of the form $a \cdot x \cdot b$ are
dense in $X_p$, so it is sensible just to say that $f \cdot x = x \cdot f$
for all $x \in X$ and $f \in C(\II)$.

For $t \in \II$, the set $I_t \coloneqq C_0(\II\setminus\{t\}, A)$ is an
ideal of $C(\II, A)$ with quotient canonically isomorphic to $A$, and each
$X_p \cdot C_0(\II \setminus \{t\}) = X_p \cdot I_t$ is a closed submodule of
$X_p$. We write $X^t_p \coloneqq X_p/(X_p\cdot I_t)$. By (the argument of)
\cite[Lemma~3.1]{FGS}, the set $\bigsqcup_p X^t_p$ is a regular product
system over $\NN^\Lambda$ of $C^*$-correspondences over $A$. Theorem~3.7 of
\cite{FGS} says that when $\Lambda$ is finite, $K_*(\Oo_{X^0}) \cong
K_*(\Oo_{X^1})$. We recover and extend this result using
Theorem~\ref{thm:special lift K-isos}.

\begin{cor}\label{cor:prodsys}
Let $A$ be a $C^*$-algebra, let $\Lambda$ be a set, and let $X$ be a homotopy
of regular product systems over $\NN^\Lambda$ of $C^*$-correspondences over
$A$. Then $K_*(\Oo_X) \cong K_*(\Oo_{X^t})$ for all $t \in \II$.
\end{cor}
\begin{proof}
Fix $t \in \II$. The argument of Proposition~3.5 of \cite{FGS} shows that
there is a surjective homomorphism $\pi : \Oo_X \to \Oo_{X^t}$ that sends
$j_{C(\II, A)}(f)$ to $j_A(f(t))$ and $j_{X}(x) = j_{X^t}(x + X_p\cdot I_t)$
for $x \in X_p$. The gauge actions on $\Oo_X$ and $\Oo_{X^t}$ induce gradings
with respect to which $\pi$ is a graded homomorphism. These are strong
gradings because $X$ is nondegenerate and full. That $X$ is nondegenerate and
full also ensures that $C \coloneqq \Oo_X$, the subalgebra $\Cc \coloneqq
j_{C(\II, A)}(C(\II, A))$ and the subspaces $j_X(X_{e_i}) \subseteq
(\Oo_X)_{e_i}$ for $1 \le i \le k$ satisfy (1)--(3) of
Theorem~\ref{thm:special lift K-isos}. Since $\pi|_\Cc$ is the evaluation map
$f \mapsto f(t)$, and since $\II$ is contractable, $\pi|_\Cc$ induces
isomorphisms $(\pi|_\Cc)_* : K_*(C(\II, A)) \to K_*(A)$. So the result
follows from Theorem~\ref{thm:special lift K-isos}.
\end{proof}

\end{document}